\documentclass{amsart}
\title{The Newton Procedure  for several variables}
\author{M. J. Soto}
\address{Departamento de \'Algebra\\Universidad de Sevilla}
\email{soto@us.es}
\author{J. L. Vicente}
\address{Departamento de \'Algebra\\Universidad de Sevilla}
\email{jlvc@us.es}

\usepackage{bm}
\theoremstyle{plain}
\newtheorem{teorema}             {Theorem}
\newtheorem{proposicion}[teorema]{Proposition}
\newtheorem{corolario}  [teorema]{Corollary}
\newtheorem{lema}       [teorema]{Lemma}
\theoremstyle{definition}
\newtheorem{definicion} [teorema]{Definition}
\theoremstyle{remark}
\newtheorem{notaciones} [teorema]{Notations}
\newtheorem{nota}       [teorema]{Remark}
\newtheorem*{inductionAssumption}{Induction assumption}
\newcommand{\vg}[2]{(#1_{1},\dots,#1_{#2})}
\newcommand{\ZZ}{\mathbb{Z}}
\newcommand{\RR}{\mathbb{R}}

\newcommand{\QQ}{\mathbb{Q}}
\newcommand{\llbrack}{\llbracket}
\newcommand{\rrbrack}{\rrbracket}
\usepackage[only,llbracket,rrbracket]{stmaryrd}
\newcommand\milex{\leq_{\text{lex}}}
\newcommand\pe[1]{\langle#1\rangle}
\newcommand\sr[1]{\pe{#1}_{+}}

\renewcommand{\qedsymbol}{\vrule height 8pt depth 1pt width 2pt}
\renewcommand{\qed}{\ifmmode\mathqed\else\unskip\nobreak\quad\hbox{\qedsymbol}\fi}
\def\lacopmbus#1#2{li\-near au\-to\-mor\-phism $#1$ of $#2$,
   which is a com\-po\-si\-ti\-on of a fi\-ni\-te
se\-quen\-ce of
   or\-der-pre\-ser\-ving mo\-no\-mial blo\-wing-ups}

\newcommand{\cS}{\mathcal{S}}
\newcommand{\RRmn}{\RR_{\geq}^{n}}

\begin{document}
\begin{abstract}
Let us consider an equation of the form
\[
P(\bm{x},z)=
z^{m}+w_{1}(\bm{x})z^{m-1}+\cdots+w_{m-1}
(\bm{x})z+w_{m}(\bm{x})=0\, ,
\]
where $m>1$, $\bm{x}=\vg{x}{n}$, $n>1$,  is a vector  of
variables, $k$ is an algebraically closed field of
characteristic zero,
$w_{i}(\bm{x})\in
k\llbrack\bm{x}\rrbrack $, and
$w_{m}(\bm{x})\neq 0$.
The aim is
to prove the Theorem of Newton-Puiseux, namely:

\begin{teorema}\label{1702061} The roots of the above
equation are formal power series with rational
exponents of bounded denominators, whose Newton diagrams
are contained in an $S$-cone.
\end{teorema}

As an application, in Section~\ref{s:closure},
we deal with some topics of integral dependence of Puiseux
power series. In particular, we
construct a domain $k\llbrack \bm{x}\rrbrack^{*}$
containing $k\llbrack \bm{x}\rrbrack$, integrally closed
in its quotioent field $k((\bm{x}))^{*}$ and this one is the
algebraic closure of~$k((\bm{x}))$.
\end{abstract}
\maketitle

\section{Introduction}{\label{s:intro}} In the last decade,
there have been several (successful) attempts to solve an
equation of integral dependence
\[
P(\bm{x},z)=
z^{m}+w_{1}(\bm{x})z^{m-1}+\cdots+w_{m-1}
(\bm{x})z+w_{m}(\bm{x})=0\, ,
\]
where $m>1$, $\bm{x}=\vg{x}{n}$, $n>1$,  is a vector  of
variables, $k$ is an algebraically closed field of
characteristic zero,
$w_{i}(\bm{x})\in
k\llbrack\bm{x}\rrbrack $, and
$w_{m}(\bm{x})\neq 0$. We will assume, in addition, that
the equation has only
simple roots in any algebraic closure of $k((\bm{x}))$,
although this means no restriction. In our opinion, the
the firt remarkable attempt is the one by
McDonald (c.f. \cite{McD}).
Later, Gonz\'{a}lez P\'{e}rez greatly extended
MacDonald's
results applying them to quasi-ordinary Puiseux power
series (c.f. \cite{GP}).
In all the cases we know, the production of the
roots is the result of a non-easy combinatorial
procedure based upon the Newton polyhedron of the whole
equation.

We have taken completely different point of view, the
simplest possible we could think. We
single out a variable, say $x_{1}$ and solve the equation in
$(z,x_{1})$ over the field $k((x_{2},\ldots,x_{n}))$ using
the elementary Newton procedure for two variables (cf.
\cite{Walker}, chapter 4, \S 3). The gain
in simplicity is enormous. The possible loss in generality
is not so much, because an usual technique in geometry is to
prepare the equations before solving them. Moreover, this
simplicity makes our techniques suitable for applications in
fields of Mathematics other than Algebra, since the tools we
use belong to the common ground of the mathematical
knowledge.

The key part  of our work is to control where the monomials
with negative exponents of the solutions lie. Surprisingly
enough,
the
Jung-Abhyankar theorem (c.f. \cite{Ab1}) gives us the clue.
In fact, in
\cite{SV2}, theorem 13, we already proved theorem
\ref{1702061},
based on the Jung-Abhyankar theorem (cited J-A from now on).
Once we know how
to control the monomials with negative exponents, we produce
here a direct proof, i.e. a proof based only on a
detailed analysis of the Newton procedure for two
variables, without ressource to J-A.

This approach is the key step to give an elementary proof of
J-A. In fact, we conjecture that one can do such a thing
 by elementary methods,
based upon theorem \ref{1702061}. We
will not deal here with such matter.

In section \ref{pps} we give a very simple description of
the roots. In section \ref{npro} we prove theorem \ref{1702061}
in the way we said above. In section
\ref{s:closure} we give the applications.

\section{Puiseux power series and $S$-cones}{\label{pps}}

 In this section we introduce
a special kind of Puiseux power series, which
will  be the roots of the
equation $P(\bm{x},z)=0$.
In other words, we are going to give meaning to the
statement of Theorem~\ref{1702061}. Let us fix an
algebraically
closed field $k$ of characteristic zero, a vector of
variables $\bm{x}=\vg{x}{n}$,
$n>1$, and an integer~$d>0$. We will consistenly use
the lexicographic order $\milex$ on $\RR^{n}$,
and the corresponding group order
 on the group of monomials
$M=\{\bm{x}^{\bm{a}}\}_{\bm{a}\in\RR^{n}}$
.

\begin{notaciones}{\label{2708051}}
Let
$\mathcal{F}_{n,d}$ be the set of all the functions
$f\colon\frac{1}{d}\ZZ^{n}\to k$;
then $\mathcal{F}_{n,d}$ is an abelian
group with respect to the usual addition of functions.
Let us write every $f\in
\mathcal{F}_{n,d}$ as a formal sum
$f=\sum_{\bm{a}\in\ZZ^{n}}f_{\bm{a}}\bm{x}^{\bm{a}/d}$
where $f_{\bm{a}}=f(\bm{a}/d)\in k$ and, if
$\bm{a}=\vg{a}{n}$, then
$\bm{x}^{\bm{a}/d}=x_{1}^{a_{1}/d}\cdots
x_{n}^{a_{n}/d}$.
We call the {\em Newton diagram} of $f$ the set
\[
\mathcal{E}(f)=
\biggl\{\frac{\bm{a}}{d}\in\frac{1}{d}\ZZ^{n}
\quad
\biggm|
\quad
\bm{a}\in\ZZ^{n} , \quad f_{\bm{a}}\neq
\bm{0}\biggr\}
.
\]

Finally, let us denote by $K_{n,d}$
the subfield $K_{n,d}=k((x_{n}^{1/d}))\dots((x_{1}^{1/d}))$
of $\mathcal{F}_{n,d}$,
which is constructed by induction. If $n=1$,  then
$K_{1,d}=k((x_{1}^{1/d}))$, the field of formal meromorphic
functions in the variable $x_{1}^{1/d}$. Any element
$\sum_{i\geq r}\alpha_{i}x_{1}^{i/d}\in k((x_{1}))$,
$\alpha_{i}\in
k$, gives the function $f\colon\frac{1}{d}\ZZ\to k$ defined
by
$f(i/d)=0$ if $i<r$ and $f(i/d)=\alpha_{i}$ for $i\geq r$.
Let us assume that $n>1$ and that we
have defined the subfield
$L=k((x_{n}^{1/d}))\cdots((x_{2}^{1/d}))$ of
$\mathcal{F}_{n-1,d}$; for each $\alpha\in L$, we denote by
$f_{\alpha}\colon\frac{1}{d}\ZZ^{n-1}\to k$ the
corresponding function.
In this situation, $K_{n,d}$ is the field
$L((x_{1}^{1/d}))$. Any element
$\sum_{i\geq r}\alpha_{i}x_{1}^{i/d}\in L((x_{1}))$,
$\alpha_{i}\in
L$, gives the function $f\colon\frac{1}{d}\ZZ^{n}\to k$
defined by
$f(i/d,a_{2}/d,\ldots,a_{n}/d)=0$ if $i<r$ and
$f(i/d,a_{2}/d,\ldots,a_{n}/d)=f_{\alpha_{i}}$ for
$i\geq r$.
\end{notaciones}

\begin{proposicion}{\label{2708052}} Let $0\neq f\in
\mathcal{F}_{n,d}$; then $f\in K_{n,d}$ if and only if\/
$\mathcal{E}(f)$ is a well-ordered subset of $\frac{1}{d}\ZZ^{n}$.
\end{proposicion}

\begin{demostracion}
Let us assume that $f\in K_{n,d}$ and use
induction on $n$.
If $n=1$, then $f\in k((x_{1}^{1/d}))$
and $\mathcal{E}(f)\subset\frac{1}{d}\ZZ$ is clearly
well-ordered. Let us assume that $n>1$ and the result
true for  $n-1$. Let
$\emptyset\neq\Omega\subset\mathcal{E}(f)$; since $f$  is a
power series in $x_{1}^{1/d}$, the set of the
first components of the vectors in $\Omega$ must have a
minimum
$a_{1}/d$. Let $0\neq u_{1}\in k((x_{n}^{1/d}))\cdots
((x_{2}^{1/d}))$ be the coefficient of $x_{1}^{a_{1}/d}$ in
$f$ and let us denote by $E$ the subset of
$\frac{1}{d}\ZZ^{n}$
consisting of all the vectors of $\mathcal{E}(u_{1})$ with
an added
$a_{1}/d$ at the beginning, as their first coordinate.
 By
the induction assumption,
$\emptyset\neq E\cap\Omega$
must have a minimum $(a_{1}/d,a_{2}/d,\ldots,a_{n}/d)$,
which is
the minimum of
$\Omega$, so $\mathcal{E}(f)$ is well-ordered.

Now, let us assume that $\mathcal{E}(f)$ is
well-ordered and use again induction on $n$. If $n=1$,
then $\mathcal{E}(f)$ has a lower bound in $\frac{1}{d}\ZZ$,
so $f\in k((x_{1}^{1/d}))$.
Let us assume that $n>1$ and the result
true for  $n-1$. Let  $a_{1}/d$ be  the first
component of the minimum of $\mathcal{E}(f)$. For a fixed
$i\in\ZZ$, $i\geq a_{1}$, we define
$u_{i}\colon\frac{1}{d}\ZZ^{n-1}\to k$ by the relation
$u_{i}(b_{2}/d,\ldots,b_{n}/d)=f(i/d,b_{2}/d,\ldots,b_{n}/d)$.
For any such $i$, the Newton diagram $\mathcal{E}(u_{i})$ is
either empty or well-ordered, so
$u_{i}\in k((x_{n}^{1/d}))\cdots((x_{2}^{1/d}))$. Therefore,
$f$ can
be written as $f=\sum_{i\geq a_{1}}u_{i}x_{1}^{i/d}$, which
implies
$f\in k((x_{n}^{1/d}))\cdots((x_{2}^{1/d}))((x_{1}^{1/d}))$.
\end{demostracion}

\begin{definicion}{\label{0604062}}
A monomial blowing-up
is a $\RR$-linear automorphism
$\varphi_{ij}$ of $\RR^{n}$,
$1\leq i\,
, j\leq n$,
$i\neq j$ defined by
\[\varphi_{ij}\vg{a}{n}=(a_{1},\ldots,\stackrel{j)}{a_{i}+a_{j}},\ldots, a_{n})\, .
\]
A monomial blowing-down is the inverse automorphism of a monomial blowing-up.
\end{definicion}

\begin{notaciones}\label{0604053}\quad
\begin{enumerate}
\item\label{0604053a0}
The {\em product order} $\ll$ is defined by
$\vg{a}{n}\ll\vg{b}{n}$ if and only if $a_{i}\leq b_{i}$,
for all $i=1,\ldots,n$.
Note that {\em the product order is
preserved by any monomial blowing-up}.

\item\label{0604053a1}
We
will also consider $\varphi_{ij}$, or its inverse, as an
automorphism of the multiplicative group $M$ of the
monomials
$\bm{x}^{\bm{a}}$,
$\bm{a}\in\RR$,  sending $\bm{x}^{\bm{a}}$ onto
$\bm{x}^{\varphi_{ij}(\bm{a})}$. This can be
viewed as applying the substitutions $x_{i}\to x_{i}x_{j}$,
$x_{l}\to x_{l}$, $l\neq i$. To apply $\varphi_{ij}$ to
$P(\bm{x},z)$
means to apply it to all its monomials, leaving $z$
fixed.

\item\label{0604053a2}
A monomial blowing-up
$\varphi_{ij}$
 preserves the lexicographic order if and only  if~$i<j$
(c.f. \cite{SV2}, proposition 5). If
$\varphi_{ij}$ is order-preserving, then so is~$\varphi_{ij}^{-1}$.
We will call them
{\em order-preserving monomial blowing-ups} or {\em order-preserving monomial blowing-downs}.
\end{enumerate}
\end{notaciones}

\begin{corolario}{\label{2708054}} Any order-preserving monomial blowing-up
$\varphi_{ij}$, $i<j$, induces a field $k$-automorphism
of $K_{n,d}$.
\end{corolario}

We borrow from \cite{SV2} (lemma 15) the
following

\begin{lema}{\label{1208052}}
Let
$\emptyset\neq\Lambda\subset\ZZ_{\geq}^{n}$;
then there exists a \lacopmbus{\Phi}{\RR^{n}},
and
a vector with integer
coordinates $\bm{a}\in\Phi(\Lambda)$ such that
$\Phi(\Lambda)\subset \bm{a}+\ZZ_{\geq}^{n}$
\end{lema}

From this result we deduce an important consequence, namely

\begin{corolario}{\label{2808052}}
Let $\Lambda_{1},\ldots,\Lambda_{r}$ be a finite number of
non-empty subsets of $\ZZ_{\geq}^{n}$. Then
there exists
a \lacopmbus{\Phi}{\RR^{n}},
and vectors
$\bm{a}_{i}\in\Phi(\Lambda_{i})$, $i=1,\ldots,r$,
such that
$\Phi(\Lambda_{i})\subset \bm{a}_{i}+\ZZ_{\geq}^{n}$.
\end{corolario}

\begin{demostracion} Let us observe that, for every
monomial blowing-up $\varphi$ and any vector
$\bm{b}\in\ZZ_{\geq}^{n}$,
one has
$\varphi(\bm{b}+\ZZ_{\geq}^{n})\subset\varphi(\bm{b})+\ZZ_{\geq}^{n}$.
Let us prove the corollary by induction on $r$. If $r=1$,
this is lemma \ref{1208052}, so let us assume that
$r>1$  and the result true for $r-1$. There exist $\Phi'$
and
$\bm{b}_{i}\in\Phi'(\Lambda_{i})$, $i=1,\ldots,r-1$,
such that
$\Phi'(\Lambda_{i})\subset\bm{b}_{i}+\ZZ_{\geq}^{n}$. On
the other hand, $\Phi'(\Lambda_{r})\subset\ZZ_{\geq}^{n}$
so, by lemma
\ref{1208052}, there exist
$\Phi''$
and $\bm{a}_{r}\in\Phi''\Phi'(\Lambda_{r})$ such that
$\Phi''\Phi'(\Lambda_{r})\subset\bm{a}_{r}+\ZZ_{\geq}^{n}$.
If, for every $i=1,\ldots,r-1$, we write
$\bm{a}_{i}=\Phi''(\bm{b}_{i})$, $i=1,\ldots,r-1$,
then, by the first observation,
\[
\Phi''\Phi'(\Lambda_{i})\subset
\Phi''(\bm{b}_{i}+\ZZ_{\geq}^{n})\subset\Phi''(\bm{b_{i}})+\ZZ_{\geq}^{n}
= \bm{a}_{i}+\ZZ_{\geq}^{n}\, .
\]
If we set $\Phi=\Phi''\Phi'$, our result is proven.
\end{demostracion}

We introduce now the objects we are looking for, namely, the
Puiseux power series in some $K_{n,d}$ whose Newton diagram is contained in an
$\mathcal{S}$-cone. In a recent paper of ours
(c.f.~\cite{SV2}), we
dealt with a special case of polyhedral cones
(see, for instance, \cite{Ewald}, page 6),
that will give rise to the $S$-cones here.

\begin{definicion}{\label{0604061}}
A
{\em polyhedral cone}
$\Gamma(\Delta)$ will be a subset of $\RR^{n}$ defined as
the projection, from the origin
$\bm{0}$, of a  compact polyhedron $\Delta$ contained in
an affine hyperplane $H$, such that $\bm{0}\notin
H$, and
$\Delta$  has a non-empty interior in $H$. In other words,
$
\Gamma(\Delta)=\cup_{\bm{a}\in\Delta}\pe{\bm{a}}_{+}
$,
where $\sr{\bm{a}}$ is the half-line of the non-negative
multiples of $\bm{a}$.
\end{definicion}

It is easy to see that the transform of  a polyhedral cone  by a monomial blowing-up,
or a monomial blowing-down, is again a polyhedral cone.
In \cite{SV2}, Theorem~6, we proved that
a polyhedral cone $\Gamma(\Delta)$ can be brought to the first quadrant
by a finite sequence of monomial blowing-ups
(i.e., its transform is contained in
$\RR_{\geq}^{n}$)
if and only if
$\Gamma(\Delta)\cap (-\RR_{\geq})^{n}=\{\bm{0}\}$.

Now we need to say more on polyhedral cones that can
be brought to the first quadrant by a finite sequence of monomial blowing-ups, namely

\begin{teorema}{\label{2008051}}
Let $\Gamma(\Delta)$ be a
polyhedral cone;
the following conditions are equivalent:
\begin{enumerate}

\item
$\Gamma(\Delta)$
can be brought to the first quadrant by a finite sequence
of order-preserving monomial blowing-ups.

\item
For every vector
$\bm{0}\neq
\bm{c}\in\Gamma(\Delta)$,
its first non-zero component is positive.

\end{enumerate}
\end{teorema}

\begin{demostracion} Let us observe that the first
non-zero component of any vector is invariant by any
order-preserving monomial blowing-up.
Consequently, if there exists a
vector $\bm{0}\neq\bm{c}\in\Gamma(\Delta)$ whose
first non-zero component is negative,
1) cannot hold.

Conversely, let us asume that 2) holds
and let
$\{\bm{c}_{1},\ldots,\bm{c}_{m}\}$, $m\geq n$, be
non-zero vectors such that the half-lines
$\sr{\bm{c}_{i}}$
 are the edges of $\Gamma(\Delta)$. Then
there must exist a finite sequence of
order-preserving monomial blowing-ups
(call
$\Phi$ their composition)  such that
$\Phi(\bm{c}_{i})\in\RR_{\geq}^{n}$, so
$\Phi(\Gamma(\Delta))\subset\RR_{\geq}^{n}$.
\end{demostracion}

\begin{definicion}{\label{2008053}} An $S$-cone is a polyhedral cone
that can be brought to the first quadrant by a finite sequence of
order-preserving monomial blowing-ups.
\end{definicion}

\begin{corolario}{\label{1808052}} Let
$\Gamma(\Delta)\subset\RR^{n}$ be an $S$-cone
and $0\neq f\in\mathcal{F}_{n,d}$ be such that
$\mathcal{E}(f)\subset\Gamma(\Delta)\cap\frac{1}{d}\ZZ^{n}$
then $f\in K_{n,d}$.
\end{corolario}

\begin{demostracion} (c.f. \cite{SV2}, proof of theorem 13).
We know that $\frac{1}{d}\ZZ_{\geq}^{n}$ is
well-ordered, so it is
$\Gamma(\Delta)$, being the inverse
image of some subset of $\frac{1}{d}\ZZ^{n}$ by a
finite sequence of order-preserving monomial blowing-downs. This implies that
$\mathcal{E}(f)$ is well-ordered and the lemma.
\end{demostracion}

\section{The Newton Procedure  for several
variables}{\label{npro}} In this section, we
construct the generalization to several variables of
the classical Newton Procedure  and prove Theorem~\ref{1702061}.
Therefore, we fix the equation $P(\bm{x},z)=0$ of the
statement of this theorem.

\begin{notaciones}{\label{2608052}} Let us consider a
polynomial
\[
Q(\bm{x}^{1/r},z)=
v_{0}(\bm{x}^{1/r})z^{m}
+
v_{1}(\bm{x}^{1/r})z^{m-1}
+\cdots+
v_{m-1}(\bm{x}^{1/r})z
+
v_{m}(\bm{x}^{1/r})
\]
where  $m\, , n\, ,
r\in\ZZ_{>}$,
$m>1$, $\bm{x}=\vg{x}{n}$ is a vector of variables,
$\bm{x}^{1/r}=(x_{1}^{1/r},\ldots,x_{n}^{1/r})$
and
\[
v_{i}(\bm{x}^{1/r})\in
k((x_{n}^{1/r}))\cdots((x_{2}^{1/r}))\llbrack x_{1}^{1/r}\rrbrack \, ,
\quad
k((x_{n}^{1/r}))\cdots((x_{2}^{1/r}))=k\; \hbox{if}
\; n=1\, ,
\]
$v_{0}(\bm{x}^{1/r})v_{m}(\bm{x}^{1/r})\neq 0$.
We  denote by $\mathcal{E}_{1}\bigl(Q(\bm{x}^{1/r},z)\bigr)$ the
{\em Newton diagram} of $Q(\bm{x},z)$ as a polynomial only in
$(x_{1},z)$, that is, we plot every monomial
$x_{n}^{a_{n}/r}\cdots x_{2}^{a_{2}/r}x_{1}^{a_{1}/r}z^{b}$,
$a_{i}\in\ZZ$, $i=1,\ldots,n$, occurring in
$Q(\bm{x}^{1/r},z)$
with a non-zero coefficient, onto the point
$(a_{1}/r,b)\in\bigl(\frac{1}{r}\ZZ_{\geq}\bigr)\times\{0,1,\ldots,m\}$
\end{notaciones}

\begin{nota}{\label{2708055}} Let $Q(\bm{x}^{1/r},z)=0$
be
a polynomial as in Notations~\ref{2608052} with $n>1$; then, for every
order-preserving monomial blowing-up (or blowing-down)
$\varphi_{ij}$, $i<j$, one has that
$\varphi_{ij}\bigl(\mathcal{E}_{1}(P)\bigr)=\mathcal{E}_{1}(P)$. The
reason is that, for any monomial
$\bm{x}^{\bm{a}/r}$, the blowing-up~$\varphi_{ij}(\bm{x}^{\bm{a}/r})$
has the same
exponent of $x_{1}$ as
$\bm{x}^{\bm{a}/r}$.
\end{nota}

The proof of Theorem\ref{1702061} is achieved by induction
on the number
$n$ of
variables in the coefficients. We
make the following
induction assumption, which
holds for
$n=1$ and
$\Phi$~equal to the identity, by the classical Theorem of
Newton-Puiseux:


\begin{inductionAssumption}[\textbf{IA}]
For every Weierstra\ss{} polynomial
with $n$~variables in the coefficient ring,
\[
\Pi(\bm{x},z)=z^{\mu}+\omega_{1}(\bm{x}^{1/\varrho})z^{\mu-1}+\cdots+\omega_{\mu-1}(\bm{x}^{1/\varrho})z+\omega_{\mu}(\bm{x}^{1/\varrho})
\in k\llbrack x_{1}^{1/\varrho},\dots,x_{n}^{1/\varrho}\rrbrack [z],
\]
with $\mu>1$,
and~$\varrho\in\ZZ_{>}$, there exists
a \lacopmbus{\Phi}{\RR^{t}},
and a
positive integer~$\pi$, such that
all the roots of $\Phi\bigl(\Pi(\bm{x},z)\bigr)=0$ belong to
$k\llbrack \bm{x}^{1/\pi}\rrbrack $.
\end{inductionAssumption}

\begin{nota}
As noted above, the case~$n=1$ is the very well known classical
Newton-Puiseux Theorem. The suite requires the reader to know in some
depth the proof, say as in~\cite{Walker}.

For our purposes, it will suffice to show here a \emph{very} brief
sketch of the methods for $n=1$ to fix some ideas and notations.

Suppose then an equation of the form
\[
\Pi(x,z)=z^{\mu}+\omega_{1}(x^{1/\varrho})z^{\mu-1}
+\dots+\omega_{\mu}(x^{1/\varrho})+
\omega_{\mu}(x^{1/\varrho}),
\]
where we will put~$\omega_{0}=1$.
If $z_{0}$ is to be a root of~$\Pi(x,z)$, we can write
\[
z_{0}=\alpha_{1}x^{\gamma_{1}}+\alpha_{2}x^{\gamma_{2}}
+\alpha_{3}x^{\gamma_{3}}+\cdots,\qquad
\gamma_{1}<\gamma_{2}<\cdots,
\]
where $\gamma_{i}\in\QQ_{>}$ for all~$i$. Rewriting
$z_{0}$ as $z_{0}=x^{\gamma}(\alpha+z_{0}')$, with $\gamma=\gamma_{1}$
and $\alpha=\alpha_{1}$, and substituting back into~$\Pi(x,z)$, we
have
\begin{align}
\Pi(x,z_{0})&=
\omega_{\mu}(x^{1/\varrho})+
\omega_{\mu-1}(x^{1/\varrho})\bigl[x^{\gamma}(\alpha+z_{0}')\bigr]
+\dots+ \bigl[x^{\gamma}(\alpha+z_{0}')\bigr]^{\mu}
\nonumber
\\
&=
\omega_{\mu}(x^{1/\varrho})+\omega_{\mu-1}(x^{1/\varrho})
x^{\gamma}\alpha+\dots+\omega_{0}x^{\gamma\mu}\alpha^{\mu}
+\Sigma(x,z_{0}'),
\label{lot}
\end{align}
where~$\Sigma(x,z_{0}')$ contains all terms on~$z_{0}'$. The idea
behind the theorem is to solve for~$\gamma$ and~$\alpha$, and iterate
the construction.

Since the order of~$z_{0}'$ is~$\gamma_{2}>0$, each term
in~$\Sigma(x,z_{0}')$ has strictly greater order than some
$\omega_{\mu-r}(x^{1/\varrho})x^{r\gamma}\alpha^{r}$. Now, a
necessary condition for~$\Pi(x,z_{0}')$ to vanish is that the lowest
order terms cancel out, so there must be at least two values of~$r$
such that
\begin{equation}
\beta=
\nu_{\mu-r_{1}}+r_{1}\gamma=
\nu_{\mu-r_{2}}+r_{2}\gamma\leq
\nu_{\mu-r}+r\gamma, \qquad
\text{for $r=0,\dots,\mu$,}
\label{ord.eq}
\end{equation}
and where $\nu_{\mu-r}$ is the order
of~$\omega_{\mu-r}(x^{1/\varrho})$.

If we group the lowest order terms in Equation~\eqref{lot}, we obtain
an equation in~$\alpha$, called the \emph{characteristic equation}, of
the form
\begin{equation}
C(\alpha)=\sum_{h} \omega_{\mu-h}'\alpha^{h},
\qquad \omega'_{\mu-h}\in k,
\label{char.eq}
\end{equation}
and where $h$ runs over all terms with $\nu_{\mu-h}+h\gamma=\beta$.

We need now to find possible values for~$\gamma$, which we do by
looking at the Newton diagram of~$\Pi(x,z)$. Equation~\eqref{ord.eq}
implies that  there exists a~$\beta$ such that all points
of~$\mathcal{E}\bigl(\Pi(x,z)\bigr)$ lie on or above the
line~$u+\gamma v=\beta$ and at least two lie exactly on it.
The linear form~$L(u,v)=u+\gamma v$ is called an
\emph{admissible linear
form} for~$\mathcal{E}\bigl(\Pi(x,z)\bigr)$. Bear in mind that the
line $u+\gamma v=\beta$ might be vertical at the very first step.

The possible values of $\gamma$ are then determined by the slopes of
the Newton polygon, and once $\gamma$~is fixed, we can solve
for~$\alpha$ in Equation~\eqref{char.eq}. Once we have $\gamma$
and~$\alpha$, we can write
\[
\Pi_{1}(x,z'')=\Pi\bigl(x,x^{\gamma_{1}}\alpha_{1}+z''\bigr),
\]
and apply the previous procedure of computing the first term to
$\Pi_{1}(x,z'')$, which is also monic in~$z''$.  The proof is
completed in~\cite{Walker} by showing that
(a)~we
can always solve for~$\alpha$ in Equation~\eqref{char.eq}, (b)~after
the very first step the Newton polygon has a segment of negative slope
and that (c)~after a finite number of steps, the $\gamma_{i}$~have a
common denominator (this is expressed by saying that the
root~$z_{0}$ has bounded denominators).
\end{nota}

It is obvious that part (a)~of the proof is trivial if we start from
an algebraically closed field~$k$, but throught the induction we will
have $\omega_{\mu-h}\in k\llbrack x_{2}^{1/\varrho},\dots,
x_{n}^{1/\varrho}\rrbrack$.
It should be noted that, since we will be
applying this very procedure for the general case,
considering
$\mathcal{E}_{1}\bigl(P(\bm{x},z)\bigr)$ and following the proof
for $n=1$, the only part we have to prove is~(a).
We do it in three lemmas.

\bigskip

Now, we start with our equation $P(\bm{x},z)=0$
and apply to it the classical
Newton Procedure, taking
$x_{1}$ and $z$ as the independent and dependent variables,
respectively. We consider the case in which the
first admissible segment is vertical.

\begin{lema}{\label{2708056}} Let us assume that
$\mathcal{E}_{1}(P)$ has a point on the vertical axis
other
than $(0,m)$ and that we choose as the first admissible
segment the vertical one. Then there exists a
\lacopmbus{\Phi_{1}}{\RR^{n}}, leaving invariant the first
coordinate of any vector, and a positive integer $d_{1}$
such that
the roots of
the corresponding~$C(\alpha)$
belong to
$k\llbrack x_{2}^{1/d_{1}},\ldots,x_{n}^{1/d_{1}}\rrbrack $.
\end{lema}

\begin{demostracion}
If we have chosen the vertical segment at the first step of
the Newton Procedure, the corresponding characteristic equation is $C(\alpha)=0$ ($\alpha$
is the unknown), where
\[
C(\alpha)=\alpha^{m}+\widehat{w}_{i_{1}}(x_{2},\ldots,x_{n})\alpha^{m-i_{1}}
+\cdots+\widehat{w}_{i_{s}}(x_{2},\ldots,x_{n})\alpha^{m-i_{s}}\, ,
\]
$s\geq 1$, $i_{1}<\cdots<i_{s}$
and
$\widehat{w}_{i_{j}}(x_{2},\ldots,x_{n})=w_{i_{j}}(0,x_{2},\ldots,x_{n})\neq 0$,
$j=1,\ldots,s$.
Then $C(\alpha)=0$ is an equation of integral dependence
over less than $n$ variables and, by {\bf IA}, there exist
$d_{1}\in\ZZ_{>}$ and a \lacopmbus{\Phi_{1}}{\RR^{n}},
leaving invariant the first coordinate of any vector,
such that the transforms by $\Phi_{1}$ of
all the roots of $\Phi_{1}\bigl(C(\alpha)\bigr)=0$ belong
to  $k\llbrack x_{2}^{1/d_{1}},\ldots,x_{n}^{1/d_{1}}\rrbrack $.
\end{demostracion}

\begin{lema}{\label{3008051}} Let
$P'(\bm{x}^{1/p},z')=0$  be an equation,
\[
P'(\bm{x}^{1/p},z')=
w'_{0}(\bm{x}^{1/p})(z')^{m}
+
w'_{1}(\bm{x}^{1/p})(z')^{m-1}
+\cdots+
w'_{m-1}(\bm{x}^{1/p})z'
+
w'_{m}(\bm{x}^{1/p})
\]
with $m>1$, $w'_{i}(\bm{x}^{1/p})\in k\llbrack \bm{x}^{1/p}\rrbrack $,
$p\in\ZZ_{>}\;$. Let us assume that
the first step
of the Newton Procedure, applied to $P'(\bm{x}^{1/p},z')$, uses any
admissible segment of negative slope in $\mathcal{E}_{1}(P')$. Under {\bf IA},
there exist a \lacopmbus{\Phi'}{\RR^{n}}, and a positive
integer $p_{1}$, such that
the roots of
the corresponding $C(\alpha)$
computed thorough the chosen segment, belong to
$k\llbrack \bm{x}^{1/p_{1}}\rrbrack $.
\end{lema}

\begin{demostracion}
Let us write
$\nu_{t}=\nu_{x_{1}}(w'_{t}(\bm{x}^{1/p}))$,
$t=0,1,\ldots,m$.
We have chosen an admissible
linear form $L=u+\gamma v$, where $u\, ,v$ are the
variables, with
$\gamma\in\QQ_{>}$, attaining a minimum
$\mu\in\QQ_{>}$ on $\mathcal{E}_{1}(P')$ at a finite
set of points
$(\nu_{m_{1}},m-m_{1})
,\ldots,
(\nu_{m_{s}},m-m_{s})$, with
$m-m_{1}>\cdots>m-m_{s}$.
Since
$\mu=\nu_{m_{t}}+\gamma(m-m_{t})$, $\forall
t=1,\ldots,s$, one must have
$\nu_{m_{1}}<\cdots<\nu_{m_{s}}$.

Let  us write
$\mathcal{E}(w'_{m_{t}}(\bm{x}^{1/p}))=\frac{1}{p}\Lambda_{m_{t}}$
with $\Lambda_{m_{t}}\subset \ZZ_{\geq}^{n}$,
for all $t=1,\ldots,s$.
By corollary \ref{2808052}, there exist a
\lacopmbus{\Phi''}{\RR^{n}}, and vectors
$\bm{a}_{m_{t}}=(a_{m_{t},1},\ldots,a_{m_{t},n})\in\Phi''(\Lambda_{m_{t}})$,
such that
$\Phi''(\Lambda_{m_{t}})\subset\bm{a}_{m_{t}}+\ZZ_{\geq}^{n}$,
that is, $\Phi''(w'_{m_{t}}(\bm{x}^{1/p}))=\bm{x}^{\bm{a}_{m_{t}}/p}
\widehat{w}'_{m_{t}}(\bm{x}^{1/p})$
where $\widehat{w}'_{m_{t}}(\bm{0})\neq 0$, for all
$t=1,\ldots,s$. We now take the equation
$\Phi''(P'(\bm{x}^{1/p},z'))=0$; one has
$\mathcal{E}_{1}(\Phi''(P'))=\mathcal{E}_{1}(P')$
by remark \ref{2708055}.

In this situation, we  have a
set of positive rationals
\[
\Omega_{1}=\{\nu_{m_{t'}}-\nu_{m_{t}}=(a_{m_{t'},1}-a_{m_{t},1})/p\mid
1\leq t<t'\leq   s\}
\]
and   a   set   of   rationals
\[
\Omega_{2}=\{(a_{m_{t},j}-a_{m_{t'},j})/p\mid
1\leq t<t'\leq s\, ,  j=2,\ldots,n\}\, .
\]
We see
that there exists a positive integer
$e$ such that each element of $e\Omega_{1}$ is greater
than
all the elements of $\Omega_{2}$.  In fact, it is enough  to
take the minimum  $\omega_{1}$ of $\Omega_{1}$,  the maximum
$\omega_{2}$  of  $\Omega_{2}$  and  $e\in\ZZ_{>}$ such that
$e\omega_{1}>\omega_{2}$.
For each $t\, ,t'$
such that $1\leq t<t'\leq s$ and each
$j=2,\ldots,n$ we have that
$e(a_{m_{t'},1}/p-a_{m_{t},1}/p)>a_{m_{t},j}/p-a_{m_{t'},j}/p$,
so
$a_{m_{t},j}/p+e a_{m_{t},1}/p<a_{m_{t'},j}/p+e
a_{m_{t'},1}/p$. Let
$\Phi''_{1}=\varphi_{1n}^{e}\cdots\varphi_{12}^{e}$ (which
clearly commute); then, if
$\bm{b}_{m_{t}}=\Phi''_{1}(\bm{a}_{m_{t}})=(b_{m_{t},1},\ldots,b_{m_{t},n})$,
one has
$\nu_{m_{t}}=a_{m_{t},1}/p=b_{m_{t},1}/p$ and
$\bm{b}_{m_{1}}\ll\bm{b}_{m_{2}}\ll\cdots\ll\bm{b}_{m_{s}}$.
Moreover,
$\Phi''_{1}\Phi''(w'_{m_{t}}(\bm{x}^{1/p}))=\bm{x}^{\bm{b}_{m_{t}}/p}
\Phi''_{1}(\widehat{w}'_{m_{t}}(\bm{x}^{1/p}))$ and this
last factor is a unit.

We now operate with the equation
$\Phi''_{1}\Phi''(P'(\bm{x}^{1/p},z'))=0$ (recall that
$\mathcal{E}_{1}(\Phi''_{1}\Phi''(P'(\bm{x}^{1/p},z'))=\mathcal{E}_{1}(P'(\bm{x}^{1/p},z'))$),
with the same chosen linear form $L$.
To compute the corresponding  terms of the
roots of $\Phi''_{1}\Phi''(P'(\bm{x}^{1/p},z'))=0$ we take
the change of variable
$z'=x_{1}^{\gamma}(\alpha+z'_{1})$
(where $\alpha$ is the unknown) and solve the characteristic equation
$C(\alpha)=0$ with
\[
C(\alpha)=\widehat{w}''_{m_{1}}(x_{2}^{1/p},\ldots,x_{n}^{1/p})\alpha^{m-m_{1}}
+\cdots+
\widehat{w}''_{m_{s}}(x_{2}^{1/p},\ldots,x_{n}^{1/p})\alpha^{m-m_{s}}\, ,
\]
where
$\widehat{w}''_{m_{t}}(x_{2}^{1/p},\ldots,x_{n}^{1/p})=
\Phi''_{1}\Phi''(w'_{m_{t}}(\bm{x}^{1/p}))
x_{1}^{-\nu_{m_{t}}}$ evaluated at
$x_{1}=0$, that is
$\widehat{w}''_{m_{t}}(x_{2}^{1/p},\ldots,x_{n}^{1/p})=
(\bm{x}^{\bm{b}_{m_{t}}/p}/x_{1}^{b_{m_{t},1}})\;
\Phi''_{1}(\widehat{w}'_{t}(0,x_{2}^{1/p},\ldots,x_{n}^{1/p}))$
and the last factor is a unit.
By the above arguments,
$\widehat{w}''_{m_{1}}(x_{2}^{1/p},\ldots,x_{n}^{1/p})$ divides all
the
other coefficients of the characteristic equation; dividing by it,
$C(\alpha)=0$ becomes an equation of integral
dependence of $\alpha$  over
$k\llbrack x_{2}^{1/p},\ldots,x_{n}^{1/p}\rrbrack $.

By {\bf IA},
there exist $p'\in\ZZ_{>}$ and a
\lacopmbus{\Phi''_{2}}{\RR^{n}} leaving invariant the first
coordinate of every vector
of $\RR^{n}$, such that, if $\alpha_{j}$,
$j=1,\ldots,r$ are the non-zero roots of
$C(\alpha)=0$, one has that
$\Phi''_{2}(\alpha_{j})\in
k\llbrack x_{2}^{1/p'},\ldots,x_{n}^{1/p'}\rrbrack $.
Writing
$\Phi'=\Phi''_{2}\Phi''_{1}\Phi''$ and
taking a common denominator $p_{1}$, we have the lemma.
\end{demostracion}

\begin{proposicion}{\label{2402061}} Let $i\geq 1$ be any
integer; there exists   a
\lacopmbus{\Phi_{i}}{\RR^{n}}, and a
positive integer $d_{i}$ such that, if we apply $i$
steps of the classical Newton Procedure  to the
equation
$\Phi_{i}\bigl(P(\bm{x},z)\bigr)=0$ in $(x_{1},z)$, in any
possible way, the sum of the first $i$ terms of any root we
obtain belongs to $k\llbrack \bm{x}^{1/d_{i}}\rrbrack $.
\end{proposicion}

\begin{demostracion}
We remind that, by Remark~\ref{2708055}, the evolution
through the Newton Procedure  of
the Newton diagram $\mathcal{E}_{1}(P)$ of $P(\bm{x},z)$ is the
same as the evolution of the Newton diagram
$\mathcal{E}_{1}\bigl(\Phi(P)\bigr)$, for any
\lacopmbus{\Phi}{\RR^{n}}; only the coefficients of the
characteristic equations change.

Lemmas \ref{2708056}~and~\ref{3008051} show that we can indeed solve
the characteristic equation in each step of the Newton procedure,
perhaps adding a composition of order-preserving blowing-ups
for every negative slope of the corresponding
Newton diagram: since all $\Phi_{i}$ are a composition of
order-preserving monomial blowing-ups, we have that the characteristic
equation of~$\Phi_{i}\bigl(P(\bm{x},z)\bigr)$ is
exactly~$\Phi_{i}\bigl(C(\alpha)\bigr)$.
With this in mind, and the fact that
any
order-preserving monomial blowing-up preserves the first quadrant of $\RR^{n}$, the proposition is
an obvious consequence of lemmas \ref{2708056} and
\ref{3008051}.
\end{demostracion}

\begin{lema}{\label{2402062}} Let
$P'(\bm{x},z')=0$  be an equation,
\[
P'(\bm{x},z')=
w'_{0}(\bm{x}^{1/p})(z')^{m}
+
w'_{1}(\bm{x}^{1/p})(z')^{m-1}
+\cdots+
w'_{m-1}(\bm{x}^{1/p})z'
+
w'_{m}(\bm{x}^{1/p})
\]
with $m>1$, $w'_{i}(\bm{x}^{1/p})\in k\llbrack \bm{x}^{1/p}\rrbrack $,
$i=0,1,\ldots,m$, $p\in\ZZ_{>}\;$. Let us assume that
$w'_{i}(0,x_{2}^{1/p},\ldots,x_{n}^{1/p})=0$, for all $i\in\{0,\ldots
m-2,m\}$ and that
$\beta=w'_{m-1}(0,x_{2}^{1/p},\ldots,x_{n}^{1/p})\neq 0$.
Then
$P'(\bm{x},z')=0$  has only one root
with positive $x_{1}$-order and
there exist a \lacopmbus{\Phi'}{\RR^{n}},
such that the transform of this root by $\Phi'$
belongs to
$k\llbrack \bm{x}^{1/p}\rrbrack $.
\end{lema}

\begin{demostracion}
The only admissible segment of $\mathcal{E}_{1}(P')$
with negative slope
consists just of the two points $(0,1)$ and
$\bigl(\nu_{x_{1}}\bigl(w'_{m}(\bm{x}^{1/p})\bigr),0\bigr)$, so the admissible
linear form is
$L=u+\gamma v$ with
$\gamma=\nu_{x_{1}}\bigl(w'_{m}(\bm{x}^{1/p})\bigr)$, and the minimum it
attains on $\mathcal{E}_{1}(P')$ is $\gamma$.
The
characteristic equation is
$0=C(\alpha)=\beta\alpha+\alpha'$,
where
$\alpha'$ is the result of making $x_{1}=0$ in
$w'_{m}(\bm{x}^{1/p})/x_{1}^{\nu_{x_{1}}(w'_{m}(\bm{x}^{1/p}))}$,
so
$\alpha=-\alpha'/\beta$.
This yields $\alpha x_{1}^{\gamma}$ as the only possible
first term of any root of $P'(\bm{x},z')=0$ with
positive $x_{1}$-order.

Now, we must perform the change of variables
$z'=x_{1}^{\gamma}(\alpha+z'_{1})$ and divide the
result by $x_{1}^{\gamma}$.
The transform of the
monomial $\beta z'$ is
$\beta \alpha x_{1}^{\gamma}+\beta x_{1}^{\gamma}z'_{1}$.
The first summand of this expression cancels with the
initial form in $x_{1}$ of $w'_{m}(\bm{x}^{1/p})$.
After this cancellation and
division by
$x_{1}^{\gamma}$ it remains the monomial
$\beta z'_{1}$,
which cannot be cancelled with any other
coming from $x_{1}^{a}(z')^{b}$
because all
of them contain a power of $x_{1}$ with exponent of the
form
$L(a,b)>\gamma$.
This shows that the transform equation is of the same form
as $P'(\bm{x},z')=0$, with the same $\beta$. This
implies
the uniqueness of the root with positive $x_{1}$-order.

Let us write
$\mathcal{E}(\beta)=\frac{1}{p}\Lambda$, with
$\emptyset\neq\Lambda\subset\ZZ_{\geq}^{n-1}$.
By lemma \ref{1208052},
there exists a \lacopmbus{\Phi'}{\RR^{n-1}},
and
a vector with integer
coordinates
$(a_{2},\ldots,a_{n})\in\Phi'(\Lambda)$ such that
$\Phi'(\Lambda)\subset
(a_{2},\ldots,a_{n})+\ZZ_{\geq}^{n-1}$, that is
$\Phi'(\beta)=
x_{2}^{a_{2}/p}\cdots x_{n}^{a_{n}/p}\beta'$,
where $\beta'$ is a unit in
$k\llbrack x_{2}^{1/p},\ldots,x_{n}^{1/p}\rrbrack $. Since all the
monomials occurring in $P'(\bm{x},z')$ contain $x_{1}$
raised to a power of the form $a/p$, $a\in\ZZ_{>}$, except
those in $\beta z'$, the same happens with
$\Phi'(P'(\bm{x},z'))$ and the exception is
$x_{2}^{a_{2}/p}\cdots x_{n}^{a_{n}/p}\beta'$.
Let $\Phi''=\varphi_{1n}^{a_{n}}\cdots\varphi_{12}^{a_{2}}$;
then all the monomials in $\Phi''\Phi'(P'(\bm{x},z'))$
are divisible by
$x_{2}^{a_{2}/p}\cdots x_{n}^{a_{n}/p}$ and only those
occurring in
$x_{2}^{a_{2}/p}\cdots x_{n}^{a_{n}/p}\beta'$ are not
divisible by $x_{1}$. Applying the Newton Procedure  to
$\Phi''\Phi'\bigl(P'(\bm{x},z')\bigr)$, as we did before to
$P'(\bm{x},z')$, it is now clear that the only root with
positive $x_{1}$-order of
$\Phi''\Phi'\bigl(P'(\bm{x},z')\bigr)=0$ belongs to
$k\llbrack \bm{x}^{1/p}\rrbrack $.
\end{demostracion}

\begin{teorema}{\label{2502061}}
There exists a positive integer $d$ and a
\lacopmbus{\Phi}{\RR^{n}}, such that
all the roots of $\Phi\bigl(P(\bm{x},z)\bigr)=0$ belong to
$k\llbrack \bm{x}^{1/d}\rrbrack $.
\end{teorema}

\begin{demostracion}
For any $i\geq 1$, proposition \ref{2402061}
tells us that our theorem is true if we consider, not the
whole roots, but the truncation of them to the first $i$
terms. In fact, this proposition tells us this result only
for some roots of the equation. Taking all the
automorphisms, composing them, and taking a common
denominator, we have the result proven for all the roots
because $\RR_{0}^{n}$ is stable by any monomial blowing-up. We know that
the classical
Newton Procedure, followed with all the necessary choices to compute
all the roots of $P(\bm{x},z)=0$  arrives at a step
in which all the equations are of the type of the one in
lemma \ref{2402062}. Composing with the new order-preserving
monomial blowing-ups given
by this lemma, we have our result.
\end{demostracion}

We finally arrive to the

\begin{proof}[Proof of Theorem \ref{1702061}]
By theorem
\ref{2502061}, the roots $\zeta_{1},\ldots,\zeta_{m}$ of
$\Phi\bigl(P(\bm{x},z)\bigr)=0$ belong to
$k\llbrack \bm{x}^{1/d}\rrbrack $. Taking into account that
$\Phi\bigl(P(\bm{x},z)\bigr)=\prod_{i=1}^{m}(z-\zeta_{i})$ and the
fact
that every monomial blowing-down is a field $k$-automorphism of $K_{n,d}$
(c.f. corollary \ref{2708054}), we have that
\[
P(\bm{x},z)=\Phi^{-1}\Phi(P(\bm{x},z))=
\prod_{i=1}^{m} \bigl(z-\Phi^{-1}(\zeta_{i})\bigr) ,
\]
so the roots of $P(\bm{x},z)$ are the $\Phi^{-1}(\zeta_{i})$, $i=1,\ldots,m$
and
$\mathcal{E}\bigl(\Phi^{-1}(\zeta_{i})\bigr)\subseteq\Phi^{-1}(\RR^{n}_{\geq})$,
which is an $\mathcal{S}$-cone.
\end{proof}

\section{Applications: integral and algebraic closures}
\label{s:closure}

Throughout this section, we will denote by $\mathcal{S}$ the
set of the finite compositions of order-preserving monomial
blowing-downs of $\RR^{n}$ and define
$
\Lambda=\{ \Phi(\RR_{\geq}^{n})\mid
\Phi\in \mathcal{S}
\}.
$
To shorten the sentences, we will simply say ``blowing-up''
(res. ``blowing-down'') instead of
order-preserving monomial
blowing-ups (resp. blowing-downs).

\begin{nota}{\label{0812062}}
Let $1\leq i<j\leq n$ be two indices and $\varphi_{ij}$
(resp. $\varphi_{ij}^{-1}$)
be the corresponding
blowing-up (resp. blowing-down); we write
$\varphi_{ij}$
(resp. $\varphi_{ij}^{-1}$)
in matrix form as $\bm{a}\to\bm{a}B$.
Then $B=E_{ij}(1)$ (resp. $B=E_{ij}(-1)$),
the elementary matrix
 equal to the $n\times n$ identity matrix $I_{n}$ except
for the fact that it has a $1$ (resp. $-1$) at the $(i,j)$
position. Then the matrix of the composition of a finite
sequence
of blowing-ups (resp. blowing-downs) has always
$1$'s at the main diagonal and it is upper-triangular.

A $n\times n$ matrix $A$
with $1$'s at the main diagonal and upper-triangular
is
 the matrix of the composition of a finite
sequence of  blowing-ups
if and only if it has
non-negative integer entries. In fact, the condition is
obviously necessary. If $A$ has non-negative integer
entries, then a suitable right-multiplication by a finite
number of matrices of the form $E_{ij}(-1)$, $i<j$, gives
the identity matrix $I_{n}$, so $A$
is
 the matrix of the composition of a finite
sequence of blowing-ups.

In the case of blowing-downs,
we can say nothing about entries. It is clear that, if a
matrix $A$ with $1$'s at the main diagonal has non-positive
entries outside it, then $A$ is the matrix
of the composition of a finite
sequence of blowing-downs, for it
can be right-multiplied by a finite sequence of matrices
$E_{ij}(1)$, $i<j$, to obtain $I_{n}$. However,
this condition being sufficient, it is not necessary:
\begin{equation}{\label{0812063}}
\left(\begin{array}{rrr}
1&-6&-8\\
0&1&-5\\
0&0&1
\end{array}\right)
\left(\begin{array}{rrr}
1&-3&-6\\
0&1&-7 \\
0&0&1
\end{array}\right)
=
\left(\begin{array}{rrr}
1&-9&28\\
0&1&-12 \\
0&0&1
\end{array}\right)\, .
\end{equation}

On the other hand, let $\bm{e}_{1},\ldots,\bm{e}_{n}$ be the
canonical base of $\RR_{\geq}^{n}$ and let
$\Phi$ be a composition of a finite sequence of
blowig-ups (resp. blowing-downs); then
\[
\Phi(\RR^{n}_{\geq})=\left\{\left.\sum_{i=1}^{n}\lambda_{i}\Phi(\bm{e}_{i})
\right|\lambda_{i}\in \RR_{\geq}\, , \,  i=1,\ldots,n
\right\}\, ,
\]
the set of non-negative linear combinations of
$\{\Phi(\bm{e}_{1}),\ldots,\Phi(\bm{e}_{n})\}$.
If $A$ is
the matrix corresponding to $\Phi$ according to the above
notations, then the row vectors of $A$ are just
$\{\Phi(\bm{e}_{1}),\ldots,\Phi(\bm{e}_{n})\}$.
\end{nota}

\begin{lema}{\label{0812061}}
For every
$\Phi\in\mathcal{S}$, all the
elements of
$\Phi(\RR_{\geq}^{n})\setminus\{\mathbf{0}\}$
are lexicographically greater than $\mathbf{0}$.
Moreover,
$\RR_{\geq}^{n}\subset\Phi(\RR_{\geq}^{n})$.
\end{lema}

\begin{proof}
The first assertion is trivial; let us show the second.
Let $A$ be the matrix of
$\Phi$; then $A^{-1}$  is the matrix of a
composition of
 blowing-ups and
$I_{n}=A^{-1}A$, which means that the vectors
of the canonical basis of
$\RR^{n}$ belong to the semigroup generated by the rows of
$A$, so
$\RR_{\geq}^{n}\subset\Phi(\RR_{\geq}^{n})$.
\end{proof}

\begin{nota}{\label{0812064}} It is not true in general
that, if $\Gamma\in\Lambda$ and $\Phi\in\mathcal{S}$ then
$\Gamma\subset\Phi(\Gamma)$
or
$\Phi(\Gamma)\subset\Gamma$. For instance, if $\Gamma$ is
given by the row vectors of the matrix $A$, $\Phi$ is given
by the matrix $B$ then $\Phi(\Gamma)$ is given by the row
vectors of the matrix $Q=AB$, where
$$
A=
\left(\begin{array}{rrr}
1&-4&-1\\
0&1&-8 \\
0&0&1
\end{array}\right)\, \quad
B=
\left(\begin{array}{rrr}
1&0&-4\\
0&1&-6\\
0&0&1
\end{array}\right)\, \quad
Q=
\left(\begin{array}{rrr}
1&-4&19\\
0&1&-14\\
0&0&1
\end{array}\right)\, ,
$$
then
$$
AQ^{-1}=
\left(\begin{array}{rrr}
1&0&-20\\
0&1&6\\
0&0&1
\end{array}\right)\, , \quad
QA^{-1}=
\left(\begin{array}{ccc}
1&0&20\\
0&1&-6\\
0&0&1
\end{array}\right)\, ,
$$
and both have a negative entry.
\end{nota}

\begin{lema}\label{star}
Let
$\Gamma_{1}=\Phi_{1}(\RR_{\geq}^{n})$,
$\Gamma_{2}=\Phi_{2}(\RR_{\geq}^{n})$,
with $\Phi_{1}$, $\Phi_{2}\in\cS$ be two $S$-cones.
Then there exists $\Phi\in\mathcal{S}$ such that
$\Gamma_{1}\subseteq\Phi(\RRmn)$ and
$\Gamma_{2}\subset\Phi(\RRmn)$.
\end{lema}

\begin{demostracion}
Let
$A_{1}$, $A_{2}$ be the respective matrices of $\Phi_{1}$
and $\Phi_{2}$; it is easy to see that
there exists
a matrix $B$, corresponding to a finite composition
of
 blowing-ups, such that $X=A_{1}B$, $Y=A_{2}B$ have their row
vectors in $\RRmn$.
The matrix
$B^{-1}$ corresponds to a $\Phi\in\mathcal{S}$ and
$A_{1}=XB^{-1}$, $A_{2}=YB^{-1}$, which means that the row
vectors of $A_{1}$, $A_{2}$ belong to the
semigroup
generated by the row vectors of $B^{-1}$. This shows that
 $\Gamma_{1}\subseteq\Phi(\RRmn)$ and
$\Gamma_{2}\subseteq\Phi(\RRmn)$.
\end{demostracion}

\begin{definicion}
If $\Gamma\in\Lambda$ and~$d\in\ZZ_{>}$,
we will write
$k\llbrack \Gamma,d\rrbrack$
for the subring of $K_{n,d}$ consisting of the Puiseux power
series whose Newton diagram is contained in $\Gamma$.
\end{definicion}

\begin{lema}  For $\Gamma\, , \Gamma'\in\Lambda$, $d\, ,
d'\in \ZZ_{>}$ one has that
\[
k\llbrack \Gamma,d\rrbrack\subset
k\llbrack \Gamma',d'\rrbrack
\iff
\text{$\Gamma\subseteq \Gamma'$ and $d|d'$}\, .
\]
Therefore, the set of rings
$k\llbrack \Gamma,d\rrbrack $,
together with the
inclusions,
is a direct system of $k$-algebras.
Also the set of their quotient fields $Q(k\llbrack
\Gamma,d\rrbrack )$,
together with the
inclusions,
is a direct system of fields.
\end{lema}

\vspace{0.3cm}

The proof is straightforward in view of lemma \ref{star}.

\begin{definicion}
The $k$-algebra $\bigcup_{(\Gamma,d)\in\Lambda\times\ZZ_{>}}
k\llbrack
\Gamma,d\rrbrack
$ will be denoted by
$k\llbrack x_{1},\dots,x_{n}\rrbrack
^{*}=k\llbrack\bm{x}\rrbrack^{*}$. If $Q(-)$ denotes
quotient fields,
 field $\bigcup_{(\Gamma,d)\in\Lambda\times\ZZ_{>}}
Q(k\llbrack \Gamma,d\rrbrack)$ will be denoted by
$k((x_{1},\dots,x_{n}))
^{*}=k((\bm{x}))^{*}$. Note that $k((\bm{x}))^{*}$ is the
quotient field of $k\llbrack \Gamma,d\rrbrack^{*}$.
\end{definicion}

We take again the Newton arguments. From theorem \ref{1702061}
we derive an easy consequence, namely the following

\begin{corolario}\label{NP2}
The roots of a polynomial
\[
P(\bm{x}^{1/d},z) =
z^{m}+\omega_{1}(\bm{x}^{1/d})z^{m-1}+\dots+\omega_{m}(\bm{x}^{1/d})\, ,
\]
where $\omega_{i}(\bm{x}^{1/d})\in
k\llbrack\bm{x}^{1/d}\rrbrack$, $i=1,\ldots,m$,
are Puiseux power series in some $K_{n,dd'}$,
such that their Newton diagrams are contained
in an~$S$-cone.
\end{corolario}

\begin{nota}
Note that Theorem~\ref{1702061} does not guarantee that all series with
exponents in~$S$-cones are
integral over~$k\llbrack \bm{x}\rrbrack $. For instance, the
power series
$
f=x^{1/2} \sqrt {1-x/y }
$
is algebraic over $k((x, y))$, its minimal polynomial is
$
z^2-x(1-x/y)\notin k\llbrack x, y\rrbrack\, ,
$
so it cannot be integral over $k\llbrack x,y\rrbrack$
because this ring is integrally closed (c.f. \cite{ZS1},
theorem 4, page 260).
\end{nota}


Finally we arrive at the results we wanted, namely

\begin{teorema}\label{integral}
The ring~$k\llbrack \bm{x}\rrbrack ^{*}$ is integrally closed (in its
quotient field).
\end{teorema}

\begin{proof}
It suffices to show that every polynomial
\[
P=z^{m}+v_{1}(\bm{x}^{1/d_{1}})z^{m-1}+\dots+v_{m}(\bm{x}^{1/d_{m}}),
\]
where~$v_{i}\in k\llbrack \Gamma_{i},d_{i}\rrbrack $,
has all its roots in~$k\llbrack \bm{x}\rrbrack ^{*}$.
Note that, by taking common denominators and common cones, we can
suppose that there exists a single pair~$(\Gamma,d)$ such
that~$v_{i}(\bm{x}^{1/d_{i}})=v_{i}(\bm{x}^{1/d})\in
k\llbrack
\Gamma,d\rrbrack
$. Now, $\Gamma=\Phi(\RRmn)$
for some
$\Phi\in\mathcal{S}$; write
\[
P'=\Phi^{-1}(P)= z^{m}+
\omega_{1}(\bm{x}^{1/d})z^{m-1}+\dots+
\omega_{m}(\bm{x}^{1/d}),
\]
with
$\omega_{i}(\bm{x}^{1/d})\in k\llbrack
\bm{x}^{1/d}\rrbrack$,
 for $i=1,\dots,m$.

We now apply Corollary~\ref{NP2} to this polynomial;
then all the roots
of~$P'$
belong to
some~$k\llbrack \Gamma',dd'\rrbrack $, where
$\Gamma'\in\Lambda$; call them
$g_{1}$, \dots,~$g_{r}$.
 But this means
that
$\Phi(g_{i})\in k\llbrack \Phi(\Gamma'),dd'\rrbrack $ are
the roots of~$P$.
\end{proof}

\begin{corolario}{\label{1012061}}
$k\llbrack\bm{x}\rrbrack^{*}$ is the integral closure of
$R=\bigcup_{\Gamma\in\Lambda}k\llbrack\Gamma,1\rrbrack$.
\end{corolario}

\begin{proof}
This proof uses well-known  facts on Galois theory of
Puiseux power series (c.f., for instance, \cite{Kiyek},
Chapter V, \S 1.2 and \S 1.3).
Any $f\in k\llbrack\Gamma,d\rrbrack$ gives raise to an
algebraic extension of $Q(R)$ generated by monomials with
exponents in $\Gamma$. By ordinary Puiseux power series
computations, the minimal polynomial of $f$ over $Q(R)$ is
and equation of integral dependence over $R$ (c.f. Kiyek,
loc.cit.). This ends our proof
\end{proof}

\begin{teorema}{\label{1012063}}
The field~$k((\bm{x}))^{*}$ is algebraically closed.
\end{teorema}

\begin{proof}
For any algebraic extension~$L=k((\bm{x}))^{*}[\alpha]$,
there exists
$c\in k\llbrack \bm{x}\rrbrack ^{*}$ such that
$\alpha'=c\alpha$ is integral over~$k\llbrack \bm{x}\rrbrack
^{*}$, and~$L=k((\bm{x}))^{*}[\alpha']$ (take a common
denominator $c\in k\llbrack \bm{x}\rrbrack ^{*}$ of the
coefficients of the minimal polynomial of $\alpha$ over
$k((\bm{x}))^{*}$). This means that
$\alpha'$~satisfies an equation of the form
\[
P=z^{m}+v_{1}(\bm{x}^{1/d_{1}})z^{m-1}+\dots+v_{m}(\bm{x}^{1/d_{m}}).
\]
By Theorem~\ref{integral}, the roots of such an equation lie
in~$k\llbrack \bm{x}\rrbrack ^{*}$, and thus, $\alpha\in k((x))^{*}$.
\end{proof}

\begin{corolario}{\label{1012062}} $k((x))^{*}$ is the
algberaic closure of $Q(R)$.
\end{corolario}

The proof is straightforward  from theorem \ref{1012063} and
corollary \ref{1012062}.

\bigskip

The first author wishes to thank the
Department of Mathematics of the University of Paderborn
for its hospitality during the writing of part of this
paper. He is especially indebted to Prof. K. Kiyek
for the fruitful discussions held during the preparation of
the manuscript.

\bibstyle{alpha}


\begin{thebibliography}{999}

\bibitem{Ab1}
S.S. Abhyankar,
{\em On the ramification of algebraic functions}.
Amer. J. Math., {\bf 36} (1955) 575--592.


\bibitem{Ewald}
G.~Ewald.
\newblock {\em Combinatorial convexity and algebraic geometry}.
Graduate Texts in Mathematics. Vol. {\bf 168},
Springer-Verlag, 1996.


\bibitem{GP}
P.D. Gonz\'{a}lez P\'{e}rez:
\newblock {\em Singularit\'{e}s quasi-ordinaires toriques
et polyh\`{e}dre de Newton du discriminant},
Canad. J. Math. Vol. {\bf 52}(2), 2000. pp. 348--368


\bibitem{Kiyek}
K. Kiyek; J.L. Vicente:
\newblock {\em Resolution of curve and surface
singularities in Characteristc Zero}, Vol {\bf 4}. Algebras
and Applications. Kluwer Academic Publishers. ISBN
1-4020-2028-7, Dordrecht, The Netherlands, 2004.

\bibitem{McD}
John McDonald,
{\em Fiber polytopes and fractional power series.}
Journal of Pure Appl. Algebra, Vol. {\bf 104},
1995 213--233.


\bibitem{SV2}
M.~J. Soto and J.~L. Vicente.
\newblock {\em Polyhedral cones and monomial blowing-ups}.
Linear Algebra and its Applications.
Vol. {\bf 412},  2006, pp. 362--372

\bibitem{Walker}
Robert J. Walker.
\newblock {\em Algebraic curves}.
\newblock Springer-Verlag.
\newblock New York.
\newblock 1978.


\bibitem{ZS1}
O.~Zariski and P.~Samuel.
\newblock {\em Commutative Algebra~I}.
\newblock Van Nostrand, 1958.


\end{thebibliography}
\end{document}